\documentclass{birkjour}
\usepackage{amssymb}
\usepackage{enumerate}

\newtheorem{thm}{Theorem}
\newtheorem{lem}[thm]{Lemma}
\newtheorem{rem}[thm]{Remark}
\newtheorem{cor}[thm]{Corollary}
\newtheorem{ex}[thm]{Example}
\begin{document}

\title{Lidskii-type formulae for Dixmier traces}
\author{Sedaev A.A.}
\address{Department of Mathematics, Voronezh State University of Architecture and Civil Engineering, 20-letiya Oktyabrya 84, Voronezh, 394006, Russia}
\email{sed@vmail.ru}
\author{Sukochev F.A.}
\address{School of Mathematics and Statistics, University of New South Wales, Kensington, NSW 2052, Australia}
\email{f.sukochev@unsw.edu.au}
\author{Zanin D.V.}
\address{School of Computer Science, Engineering and Mathematics, Flinders University,
Bedford Park, SA 5042, Australia}
\email{zani0005@csem.flinders.edu.au}

\begin{abstract} We establish several analogues of the classical Lidskii Theorem for some special classes of singular traces (Dixmier traces and Connes-Dixmier traces) used in noncommutative geometry.
\end{abstract}

\thanks{The first author was partially supported by RFBR 08-01-00226}

\maketitle

\section{Introduction and Preliminaries}

\subsection{Dixmier-Macaev ideal and Dixmier traces}

An important role in noncommutative geometry \cite{Connes} is played by the set of compact operators whose partial sums of singular values are logarithmically divergent. This set can be adequately described using the terminology of Marcinkiewicz spaces. Consider the Marcinkiewicz sequence space
$$m_{1,\infty}:=\{x=\{x_n\}_{n=1}^\infty:\|x\|_{m_{1,\infty}}<\infty\},$$
where we set
$$\|x\|_{m_{1,\infty}}=\sup_{N}\frac1{\log(N+1)}\sum_{n=1}^Nx_n^*.$$
Here, $\{x_n^*\}_{n=1}^\infty$ is the sequence $\{|x_n|\}_{n=1}^\infty$ rearranged in nonincreasing order.

Fix an infinite-dimensional separable complex Hilbert space $H$ and consider the set $\mathcal{M}_{1,\infty}$ of all compact operators $x$ on $H$ such that the sequence of its singular values $\{s_n(T)\}_{n=1}^\infty$ falls into the space $m_{1,\infty}$ (recall that the singular values of a compact operator $T$ are the eigenvalues of the operator $|T|=(T^*T)^{1/2}$). We set
$$\|T\|_{\mathcal{M}_{1,\infty}}:=\|\{s_n(T)\}\|_{m_{1,\infty}}.$$
It is well known that the ideal of compact operators $\mathcal{M}_{1,\infty}$ equipped with the norm $\|\cdot\|_{\mathcal{M}_{1,\infty}}$ is a Banach space. We refer to the recent paper \cite{Pie} by Pietsch for additional references and information on these spaces.

We describe briefly a construction of singular traces on the ideal $\mathcal{M}_{1,\infty}$ due to Dixmier \cite{Dixmier} and its various modifications which are of importance in noncommutative geometry \cite{Connes}. For a more detailed treatment we refer to \cite{CS}.

Let $\sigma_n$, $n\geq1$ be the operator on $l_{\infty}$ defined by
$$\sigma_n(x_1,\ldots,x_k,\ldots)=(\underbrace{x_1,\ldots,x_1}_{n\text{-times}},\underbrace{x_2,\ldots,x_2}_{n\text{-times}},\ldots,\underbrace{x_k,\ldots,x_k}_{n\text{-times}},\ldots).$$
Let $\omega$ be a $\sigma_n$-invariant generalised limit on $l_{\infty},$ that is, $\omega$ is a positive normalised functional on $l_{\infty}$ such that $\omega(\sigma_n(x))=\omega(x)$ for all $x\in l_{\infty}$ and such that $\omega|_{c_0}=0,$ where $c_0$ is the subspace of all vanishing sequences. For an element $0\leq T\in\mathcal{M}_{1,\infty}$ we set
$$\tau_{\omega}(T):=\omega(\{\frac1{\log(N+1)}\sum_{n=1}^N s_n(T)\}_{N=1}^\infty).$$
It is well known (see e.g. \S~5 in \cite{CS} and additional references therein) that $\tau_{\omega}$ is an additive functional on the positive part of $\mathcal{M}_{1,\infty}.$ Thus, $\tau_{\omega}$ admits a linear extension to a unitarily invariant functional (trace) on $\mathcal{M}_{1,\infty}.$ This trace vanishes on all finite-dimensional operators from $B(H).$ Such singular traces are called {\it Dixmier traces} (see \cite{Dixmier}).

A smaller subclass of Dixmier traces was introduced by Connes in \cite{Connes} by observing that a functional $\omega=\gamma\circ M$ is $\sigma_n$-invariant state on
$l_\infty$ for all $n\ge 1.$ Here, $\gamma$ is an arbitrary generalised limit on the space $L_{\infty}(0,\infty)$ of all bounded measurable functions and the operator $M$ is a Cesaro operator defined by the formula
$$(Mx)(t)=\frac1{\log(t)}\int_1^t\frac{x(s)ds}{s}.$$
Referring to $\omega$ above as a functional on $l_{\infty},$ we tacitly apply an isometric embedding $i:l_{\infty}\to L_{\infty}(0,\infty)$ given by $$\{x_j\}_{j=1}^\infty\stackrel{i}{\mapsto}\sum_{j=1}^\infty x_j \chi_{[j-1,j)},$$
where $\chi_{[j-1,j)}$ is the characteristic function of the interval $[j-1,j).$ Dixmier traces $\tau_{\omega}$ defined such $\omega$'s are termed {\it Connes-Dixmier traces}. We refer to \cite {Connes} and \cite {LSS,CS} for discussion of their properties.

Finally, various formulae of noncommutative geometry (in particular, those involving heat kernel estimates and generalised $\zeta-$function) were established in \cite{CPS,CRSS,Connes} for yet a smaller subset of Connes-Dixmier traces, when the functional $\omega$ was assumed to be $M$-invariant. This class (and its further modifications) was first introduced in \cite{CPS} (see also \cite{DPSSS2}) and further studied and used in \cite{BF,AS,CRSS}. For brevity we refer to the latter class (a proper subclass of Connes-Dixmier traces) as a class of {\it $M$-invariant Dixmier traces}.

\subsection{Lidskii formula for $M$-invariant Dixmier traces in \cite{CPS,AS,BF}}

In the case, when we deal with the standard trace ${\rm Tr}$ and the standard trace class $\mathcal {S}_1$ of compact operators from $B(H),$ the classical Lidskii Theorem asserts that the trace
$${\rm Tr}(T)=\sum_{n\geq 1}\lambda_n(T)$$
for any $T\in\mathcal{S}_1.$ Here, $\{\lambda_n(T)\}_{n\geq 1}$ is the sequence of eigenvalues of $T,$ taken in any order. This arbitrariness of the order is due to the absolute convergence of the series $\sum_{n\geq 1}|\lambda_n(T)|.$ In particular, we can choose the decreasing order of absolute values of $\lambda_n(T)$ and counting multiplicities.

The core difference of this situation with the setting of Dixmier traces living on the ideal $\mathcal{M}_{1,\infty}$ consists in the fact that the series $\sum_{n\geq 1}
|\lambda_n(T)|$ generally speaking diverges for every $T\in\mathcal{M}_{1,\infty}.$ For simplicity, we explain the emerging obstacle in the case of a self-adjoint operator $T=T^*=T_+-T_-\in\mathcal{M}_{1,\infty}.$ For such $T,$ by the definition, $\tau_{\omega}(T)=\tau_{\omega}(T_+)-\tau_{\omega}(T_-),$ where $$\tau_{\omega}(T_{\pm})=\omega(\{\frac{1}{\log(N)}\sum_{n=1}^N\lambda_n(T_\pm)\}).$$
Even in this case, it is not clear why the equality
$$\tau_\omega(T)=\omega(\{\frac{1}{\log(N)}\sum_{n=1}^N\lambda_n(T)\})$$
should hold for the special enumeration of the set $\{\lambda_n(T)\}_{n\geq 1}$ given by the decreasing order of absolute values of $|\lambda_n(T)|;$ or for that matter for \emph{any} enumeration of this set.

The following result from \cite{AS} establishes the equality above under significant additional constraints on $\tau_{\omega}$ and $T\in\mathcal{M}_{1,\infty}.$
\begin{thm}\label{formula lidskogo from AS} Let $\omega$ be $M$-invariant and let $T\in \mathcal{M}_{1,\infty}$ satisfy the assumption $s_n(T)\leq C/n$ for some $C>0$ and all
$n\ge 1$. We have
$$\tau_\omega(T)=\omega(\frac{1}{\log(n)}\sum_{|\lambda|>1/n,\lambda\in\sigma(T)}\lambda),$$
where $\sigma(T)$ is the spectrum of $T.$
\end{thm}

In the case when $T$ is a {\it positive} arbitrary element from $\mathcal{M}_{1,\infty}$ and $\omega$ is taken from a rather special subset of all $M$-invariant generalised limits (termed
in \cite{CPRS1} "maximally invariant Dixmier functionals") this result can be already found in \cite[Proposition 2.4]{CPS}. In \cite[Theorem 1]{BF}, the assertion from \cite[Proposition
2.4]{CPS} was extended to an arbitrary $M$-invariant $\omega.$ Another modification of the class of $\omega$'s for which the result of \cite[Proposition 2.4]{CPS} and \cite[Theorem 1]{BF} holds is given in \cite[Proposition 4.3]{CRSS}.

\subsection{Statement of main results}

In this paper we prove significant extensions and generalisations of Theorem \ref{formula lidskogo from AS} from \cite{AS}, \cite[Proposition 2.4]{CPS}, \cite[Theorem 1]{BF} and \cite[Proposition 4.3]{CRSS}. Many of our results are established for a general class of Marcinkiewicz ideals. Here, for convenience of the reader, we restate these results for
traces on $\mathcal{M}_{1,\infty}.$ Our first main result shows that the assertion of Theorem \ref{formula lidskogo from AS} holds for an arbitrary Connes-Dixmier trace $\tau_{\omega}$.

\begin{thm}\label{formulalidskogo1} Let $\tau_{\omega}$ be a Connes-Dixmier trace on $\mathcal{M}_{1,\infty}.$ We have
\begin{equation}\label{as lidskiy formula}
\tau_\omega(T)=\omega(\frac{1}{\log(n)}\sum_{|\lambda|>1/n,\lambda\in\sigma(T)}\lambda),\quad
T\in\mathcal{M}_{1,\infty}.
\end{equation}
\end{thm}

Theorem \ref{formulalidskogo1} follows immediately from Theorem \ref{formula lidskogo dlya l1infty} below.

Our second main result is the answer to a natural question whether formula \eqref{as lidskiy formula} holds for every Dixmier trace. This question is answered in negative in Theorem
\ref{fails for Dixmier traces}.

Our third (and the last) main result answers in the affirmative the question whether there exists a modification of the summation method used in formula \eqref{as lidskiy formula}
ensuring that it holds for all Dixmier traces.

\begin{thm}\label{formulalidskogo2} Let $\tau_{\omega}$ be a Dixmier trace on $\mathcal{M}_{1,\infty}.$ We have
$$\tau_\omega(T)=\omega(\frac1{\log(n)}\sum_{\lambda\in\sigma(T),|\lambda|>\log(n)/n}\lambda),\quad T\in \mathcal{M}_{1,\infty}.$$
\end{thm}

Theorem \ref{formulalidskogo2} follows immediately from Theorem \ref{lidskii formula} below.

At the end of the paper we also provide an application of our results. The result proved in the last section concerns heat kernel type formulae from noncommutative geometry (see \cite{Connes, CPS, CRSS, BF}) and has been already established in \cite{Sedaev} with a rather arcane argument.  We present here a very simple approach to these formulae.

\subsection{Marcinkiewicz spaces and singular traces}

It is convenient to consider the general class of Marcinkiewicz spaces since many of our results hold for this class with no extra effort. We frequently use commutative results as a stepping stone to obtain their noncommutative analogues.

Recall that the distribution function $n_x$ of a bounded measurable function $x$ is defined by the formula
$$n_x(t)=m(\{s,\ |x(s)|>t,\quad t>0\}).$$
We write $x^*$ for the decreasing rearrangement of the function $x$: $x^*$ is the right continuous non-increasing function whose distribution function coincides with that of $|x|$ (see \cite{KPS}).

The following formula is frequently used in the proofs below sometimes without explicit referencing.
\begin{equation}\label{eq3}
\int_0^{n_x(t)}x^*(s)ds=-\int_t^{\infty}\lambda dn_x(\lambda).
\end{equation}
Here, $z$ is any positive number.

Marcinkiewicz spaces are a special case of fully symmetric function  and sequence spaces, see \cite{KPS}. Denote by $\Psi$ the class of all concave increasing functions such that $\psi(\infty)=\infty,$ $\psi(t)=O(t)$ as $t\to0$ and $\psi(t)=o(t)$ as $t\to\infty.$ For every $\psi\in\Psi,$ Marcinkiewicz space $M_{\psi}$ is a set of all bounded measurable functions $x$ on $[0,\infty)$ such that
\begin{equation}
\|x\|_{M_{\psi}}:=\sup_{t>0}\frac1{\psi(t)}\int_0^tx^*(s)ds<\infty.
\end{equation}
Marcinkiewicz sequence space $m_{\psi}$ is a set of sequences (see e.g. \cite{Pie,CS}) satisfying the condition
$$\|x\|_{m_{\psi}}=\sup_n\frac1{\psi(n)}\sum_{k=1}^nx_n^*<\infty.$$

In this paper, we mainly work with functions $\psi\in\Psi$ satisfying the following condition.
\begin{equation}\label{good upper bound}
\limsup_{t\to\infty}\frac{\psi(2t)}{\psi(t)}<2.
\end{equation}

Let $K(H)$ be the ideal of all compact operators. If $m_{\psi}$ is a Marcinkiewicz sequence space, then the corresponding Marcinkiewicz operator space $\mathcal{M}_{\psi}$ is the set of all $T\in K(H)$ such that $\{s_n(T)\}\in m_{\psi}$ equipped with the norm $\|T\|_{\mathcal{M}_{\psi}}:=\|\{s_n(T)\}\|_{m_{\psi}}.$

Let $\psi\in\Psi$ and let $\omega$ be a dilation invariant generalised limit. The mapping
$\tau_{\omega}$ defined by the formula
$$\tau_{\omega}(x):=\omega(\frac1{\psi(t)}\int_0^tx^*(s)ds)$$
is a subadditive homogeneous functional on $M_{\psi}^+.$ If $\tau_{\omega}$ is additive on $M_{\psi}^+,$ then $\tau_{\omega}$ is called {\it Dixmier trace generated by $\omega$}. We refer the reader to \cite{DPSS,DPSSS,DPSSS2} for conditions which guarantee the additivity of $\tau_{\omega}.$ It is well known that $\tau_{\omega}$ is additive for any $\omega$ as above when
\begin{equation}\label{limit condition}
\lim_{t\to\infty}\frac{\psi(2t)}{\psi(t)}=1.
\end{equation}
Similarly, the definitions of Connes-Dixmier traces and $M$-invariant traces naturally extend to denote corresponding singular traces on Marcinkiewicz ideals $\mathcal{M}_{\psi}$ (see \cite{LSS}).

Our main result for general Dixmier traces on ideals
$\mathcal{M}_{\psi}$ is given in Theorem \ref{lidskii formula} which asserts that for any Dixmier trace $\tau_{\omega}$ on $\mathcal{M}_{\psi}$ with $\psi\in\Psi$ satisfying condition \eqref{good upper bound} we have
\begin{equation}\label{general lidskii formula}
\tau_\omega(T)=\omega(\frac1{\psi(n)}\sum_{\lambda\in\sigma(T),|\lambda|>\psi(n)/n}\lambda),\quad T\in \mathcal{M}_{\psi}.
\end{equation}

The result of Theorem \ref{formulalidskogo2} follows immediately from the formula above, if we set $\psi(t)=\log(t)$ for all $t\geq 2.$

\subsection{Failure of \eqref{as lidskiy formula} for Dixmier traces}

Here, we show that there are Dixmier traces $\tau_{\omega}$ on $\mathcal{M}_{1,\infty}$ for which formula \eqref{as lidskiy formula} fails. To this end we use $\omega$ provided by the lemma below.

Define a subadditive functional $\pi:L_{\infty}(0,\infty)\to\mathbb{R}$ by the formula
$$\pi(x)=\limsup_{N\to\infty}\frac1{\log(\log(N))}\int_N^{N\log(N)}\frac{x(s)ds}{s}.$$
Clearly, $\pi$ is positive and homogeneous.

The following lemma is routine. We include the proof for convenience of the reader.

\begin{lem}\label{pi lemma}
\begin{enumerate} Let $x\in L_{\infty}(0,\infty)$ be an arbitrary positive element.
\item If $\omega\in L_{\infty}(0,\infty)^*$ such that $\omega\leq\pi,$ then $\omega$ is dilation invariant generalised limit.
\item If $\pi(x)>0,$ then there exists a dilation invariant generalised limit $\omega$ such that $\omega(x)>0.$
\end{enumerate}
\end{lem}
\begin{proof} We prove the first assertion and then derive the second one from it.
\begin{enumerate}
\item At first we note that by assumption
\begin{equation}\label{up down estimate}
-\pi(-y)\leq\omega(y)\leq\pi(y)
\end{equation}
for every $y\in L_{\infty}(0,\infty).$ Note that $\pi(-y)\leq0$ for every $0\leq y\in L_{\infty}(0,\infty).$ It follows that $\omega$ is positive.

Further, for every $y\in L_{\infty}(0,\infty)$ we have
$$|\int_N^{N\log(N)}\frac{(y-\sigma_ny)(s)ds}{s}|=|\int_{N\log(N)/n}^{N\log(N)}\frac{y(s)ds}{s}-\int_{N/n}^N\frac{y(s)ds}{s}|.$$
Therefore,
$$|\pi(y-\sigma_ny)|\leq\limsup_{N\to\infty}\frac1{N\log(N)}\cdot 2\|y\|_{\infty}\cdot |\log(n)|=0.$$
Hence,
$$\omega(y-\sigma_ny)\leq\pi(y-\sigma_ny)=0,\quad \omega(-y)\leq\pi(\sigma_ny-y)=0.$$
Thus, $\omega$ is dilation invariant.

If $y\in L_{\infty}(0,\infty)$ is such that $y(t)\to0$ as $t\to\infty,$ then $\pi(y)=\pi(-y)=0.$ It follows from \eqref{up down estimate} that $\omega(y)=0.$

Noting that $\omega(1)=1,$ we conclude that $\omega$ is dilation invariant generalised limit.

\item Consider linear space $x\mathbb{R}$ spanned by element $x.$ Set $\omega(\lambda x)=\lambda\pi(x)$ for every $\lambda\in\mathbb{R}.$ It follows that $\omega\leq\pi$ on $x\mathbb{R}.$ By the Hahn-Banach theorem, there exists a functional $\omega\in L_{\infty}(0,\infty)^*$ such that $\omega(x)=\pi(x)$ and $\omega\leq\pi.$ It follows from above that $\omega$ is a dilation invariant generalised limit.
\end{enumerate}
\end{proof}

\begin{thm}\label{fails for Dixmier traces} There exist a positive function $x\in M_{1,\infty}$ and a Dixmier trace $\tau_{\omega}$ such that
\begin{equation}\label{no lidskiy}
\tau_{\omega}(x)\neq\omega(\frac{-1}{\log(t)}\int_{1/t}^{\infty}\lambda dn_x(\lambda)).
\end{equation}
\end{thm}
\begin{proof} Define a function $x$ by the formula
$$x=\sup_k e^{-e^k}\chi_{[1,e^{k+e^k}]}.$$
If $t\in[e^{k-1+e^{k-1}},e^{k+e^k}],$ then
$$\frac1{\log(t)}\int_1^tx^*(s)ds\leq e^{1-k}\int_1^{e^{k+e^k}}x^*(s)ds\leq e^{1-k}\sum_{n=1}^ke^{-e^n}\cdot e^{n+e^n}\leq\frac{e^2}{e-1}.$$
Thus, $x\in M_{1,\infty}.$

We claim that
$$\limsup_{N\to\infty}\frac1{\log(\log(N))}\int_N^{N\log(N)}(\frac1{\log(t)}\int_t^{n_x(1/t)}x^*(s)ds)dt>0.$$
Set $N=e^{e^k}.$ It is clear that $n_x(1/t)=e^{k+e^k}$ for every $t\in [N,N\log(N)].$ Since $x^*(s)=e^{-e^k}$ for every $s\in[t,n_x(1/t)]$ and every $t\in[N,N\log(N)],$ we can rewrite the expression under the limit in the left-hand side as
$$\frac1k\int_{e^{e^k}}^{e^{k+e^k}}\frac{e^{k+e^k}-t}{e^{e^k}t\log(t)}dt=\frac{e^k}{k}\int_{e^{e^k}}^{e^{k+e^k}}\frac{dt}{t\log(t)}-\frac1{ke^{e^k}}\int_{e^{e^k}}^{e^{k+e^k}}\frac{dt}{\log(t)}=$$
$$=\frac{e^k}{k}\log(1+\frac{k}{e^k})-o(1)=1+o(1).$$
This proves the claim.

Thus,
$$\pi(\frac1{\log(t)}(\int_1^{n_x(1/t)}x^*(s)ds-\int_1^tx^*(s)ds))>0.$$
The assertion of the theorem now follows from Lemma \ref{pi lemma} and \eqref{eq3}.
\end{proof}

\section{Lidskii formula for Connes-Dixmier traces}

In this section, we extend results of \cite{AS} (and, partially, those of \cite{BF}) to a wider class of Marcinkiewicz spaces and Connes-Dixmier traces. To this end, we need some extra assumptions on $\psi\in\Psi.$ The need of such additional conditions is seen from the example below, which shows that analogue of formula \eqref{as lidskiy formula} for an arbitrary $\psi\in\Psi$ fails.

\begin{ex} Let $\psi(t)=\exp(\sqrt{\log(t)})$ and let $x=\psi'.$ If $\tau_{\omega}$ is a Dixmier trace on $M_{\psi},$ then
$$e^{1/2}\tau_{\omega}(x)\leq\omega(\frac{-1}{\psi(t)}\int_{1/t}^{\infty}\lambda dn_x(\lambda)).$$
\end{ex}
\begin{proof} It is clear that $x(t)=\exp(\sqrt{\log(2)})/2t\sqrt{\log(t)}.$ We have
$$\frac{t\exp(\sqrt{\log(t)})}{2\sqrt{\log(t)}}\leq n_x(1/t).$$
for all sufficiently large $t.$ Hence,
$$e^{1/2}+o(1)\leq\frac{\psi(n_x(1/t))}{\psi(t)}.$$
The assertion follows immediately.
\end{proof}

Thus, some additional restrictions on the function $\psi$ are needed. We require the following condition
\begin{equation}\label{sedaev psi}
\lim_{t\to\infty}\frac{\psi(t\psi(t))}{\psi(t)}=1.
\end{equation}
It is clear that \eqref{limit condition} holds and, therefore, Marcinkiewicz space $M_{\psi}$ admits nonzero Dixmier traces (see \cite{DPSS},\cite{DPSSS},\cite{DPSSS2}).

Now we show that formula \eqref{as lidskiy formula} holds for all Connes-Dixmier traces on $\mathcal{M}_{\psi}.$

\begin{lem}\label{second dx estimate} Let $\psi\in\Psi$ satisfy condition \eqref{sedaev psi}. If $c>\|x\|_{M_{\psi}},$ we have
$$d_x(1/t)\leq ct\psi(t)$$
for every $x\in M_{\psi}$ and every sufficiently large $t.$
\end{lem}
\begin{proof} Assume the contrary. Hence, there exists a sequence $t_k\to\infty$ such that $x^*(s)\geq 1/t_k$ for every $s\in[0,ct_k\psi(t_k)].$ By the definition of  Marcinkiewicz norm,
$$\|x\|_{M_{\psi}}\geq\frac1{\psi(ct_k\psi(t_k))}\int_0^{ct_k\psi(t_k)}x^*(s)ds\geq\frac{c\psi(t_k)}{\psi(ct_k\psi(t_k))}.$$
It follows from \eqref{sedaev psi} that
$$\frac{c\psi(t_k)}{\psi(ct_k\psi(t_k))}\to c.$$
The contradiction proves the Lemma.
\end{proof}

\begin{rem}\label{1+o1 remark} Let $0\leq x,y\in L_{\infty}(0,\infty)$ and let $y(t)=x(t)\cdot(1+o(1))$ as $t\to\infty.$ If $x\notin L_1(0,\infty),$ we have
$$\int_1^Ty(s)ds=(1+o(1))\int_1^Tx(s)ds.$$
\end{rem}

\begin{lem}\label{M k nulyu} Let $\psi\in\Psi$ satisfy condition \eqref{sedaev psi}. We have
$$\frac1{\log(T)}\int_1^T\frac{dt}{t\psi(t)}\int_0^{ct\psi(t)}x(s)ds=\frac1{\log(T)}\int_1^T\frac{dt}{t\psi(t)}\int_0^{t}x(s)ds+o(1)$$
as $T\to\infty$ for every positive $x=x^*\in M_{\psi}$ and every $c>0.$
\end{lem}
\begin{proof} The assertion is linear with respect to $x.$ Since the assertion holds for $x(t)=\psi'(t),$ it is sufficient to verify it for $x+\psi'$ instead of $x.$ Hence, we may assume that $x(t)\geq\psi'(t).$ Thus, integral in the right-hand side is unbounded as $T\to\infty.$

Make a substitution $z=ct\psi(t)$ in the left-hand side integral. It follows from the condition \eqref{sedaev psi} that
$$\frac{dt}{t\psi(t)}=\frac{dz}{z\psi(z)}(1+o(1)).$$
Indeed, by Lagrange theorem, we have
$$\psi(z)=\psi(t)(1+o(1)),\quad \frac{dz}{z}=\frac{dt}{t}(1+\frac{t\psi'(t)}{\psi(t)})=\frac{dt}{t}(1+o(1)).$$
It follows from Remark \ref{1+o1 remark} that
\begin{equation}\label{eq1}
\int_1^T\frac{dt}{t\psi(t)}\int_0^{ct\psi(t)}x^*(s)ds=(1+o(1))\int_{c\psi(1)}^{cT\psi(T)}\frac{dz}{z\psi(z)}\int_0^zx^*(s)ds.
\end{equation}
Evidently,
\begin{equation}\label{eq2}
\int_{T}^{cT\psi(T)}\frac{dz}{z\psi(z)}\int_0^zx^*(s)ds=O(\int_{T}^{cT\psi(T)}\frac{dz}{z})=o(1).
\end{equation}
Noting that
$$\int_1^{cT\psi(T)}=\int_1^T+\int_T^{cT\psi(T)}$$
the combination of \eqref{eq1} and \eqref{eq2} yields the assertion.
\end{proof}

\begin{lem}\label{ocenka snizu connes} Let $\psi\in\Psi$ satisfy the condition \eqref{sedaev psi} and let $\tau_{\omega}$ be a Connes-Dixmier trace on $M_{\psi}.$ We have
$$\omega(\frac{-1}{\psi(t)}\int_{1/t}^{\infty}\lambda dn_x(\lambda))\leq\tau_{\omega}(x)$$
for every positive $x\in M_{1,\infty}.$
\end{lem}
\begin{proof} Due to \eqref{eq3} and Lemma \ref{second dx estimate} we have
$$\omega(\frac1{\psi(t)}\int_0^{n_x(1/t)}x^*(s)ds)\leq\omega(\frac1{\psi(t)}\int_0^{ct\psi(t)}x^*(s)ds)=$$
$$=\gamma(M(\frac1{\psi(t)}\int_0^tx^*(s)ds)+o(1))=\gamma(M(\frac1{\log(t)}\int_0^tx^*(s)ds))=\tau_{\omega}(x).$$
\end{proof}

\begin{lem}\label{prostaya ocenka sverhu} Let $\psi\in\Psi$ and let $\tau_{\omega}$ be a Dixmier trace on $M_{\psi}.$ We have
$$\tau_{\omega}(x)\leq\omega(\frac{-1}{\psi(t)}\int_{1/t}^{\infty}\lambda dn_x(\lambda))$$
for every positive $x\in M_{\psi}.$
\end{lem}
\begin{proof} We claim that
$$\int_0^tx^*(s)ds\leq\int_0^{n_x(1/t)}x^*(s)ds+1.$$
The inequality is evident if $t\leq n_x(1/t).$ If $t>n_x(1/t),$ then $x^*(s)\leq 1/t$ for every $s\in[n_x(1/t),t].$ It follows that
$$\int_0^tx^*(s)ds=\int_0^{n_x(1/t)}x^*(s)ds+\int_{n_x(1/t)}^tx^*(s)ds\leq$$
$$\leq\int_0^{n_x(1/t)}x^*(s)ds+(t-n_x(1/t))\cdot t^{-1}.$$
Thus, claim holds in either case.

It follows that
$$\tau_{\omega}(x)\leq\omega(\frac1{\psi(t)}\int_1^{n_x(1/t)}x^*(s)ds)+\omega(\frac1{\psi(t)}).$$
The assertion follows immediately.
\end{proof}

The next theorem follows immediately from Lemma \ref{ocenka snizu connes} and Lemma \ref{prostaya ocenka sverhu}.

\begin{thm}\label{formula lidskogo dlya l1infty comm} Let $\psi\in\Psi$ satisfy the condition \eqref{sedaev psi} and let $\tau_{\omega}$ be a Connes-Dixmier trace on $M_{\psi}.$ We have
$$\tau_\omega(x)=\omega(\frac{-1}{\psi(t)}\int_{1/t}^{\infty}\lambda dn_x(\lambda))$$
for every positive  $x\in M_{\psi}.$
\end{thm}

\begin{rem} Consider weak space $M_{\psi}^w$ (the smallest symmetric ideal containing $\psi'$). Suppose that $\psi$ satisfies the condition \eqref{sedaev psi}. If $\tau_{\omega}$ is an arbitrary Dixmier trace on $M_{\psi},$ then we have
$$\tau_{\omega}(x)=\omega(\frac{-1}{\psi(t)}\int_{1/t}^{\infty}\lambda dn_x(\lambda))$$
for every positive $x\in M_{\psi}^w.$ Using Lemma \ref{second dx estimate}, the equality above follows immediately.
\end{rem}

Arguing as in the section \ref{noncomm section} below, we obtain a noncommutative version of Theorem \ref{formula lidskogo dlya l1infty comm}, which strengthens \cite[Theorem 1]{BF} and \cite[Corollary 2.12]{AS} (see Theorem \ref{formula lidskogo from AS}).

\begin{thm}\label{formula lidskogo dlya l1infty} Let $\psi\in\Psi$ satisfy the condition \eqref{sedaev psi} and let $\tau_{\omega}$ be a Connes-Dixmier trace on $\mathcal{M}_{\psi}.$ We have
$$\tau_\omega(T)=\omega(\frac{1}{\psi(n)}\sum_{|\lambda|>1/n,\lambda\in\sigma(S)}\lambda)$$
for every operator $T\in\mathcal{M}_{\psi}.$
\end{thm}

\section{Adjusted Lidskii formula for Dixmier traces: commutative setting}

As we have seen in Theorem \ref{fails for Dixmier traces}, formula \eqref{as lidskiy formula} does not hold for Dixmier traces $\tau_{\omega}.$ In this section, we consider a modification of formula \eqref{as lidskiy formula} which holds for all Dixmier traces $\tau_{\omega}$ on a commutative Marcinkiewicz space $M_{\psi},$ $\psi\in\Psi.$

\begin{lem}\label{dx estimate} Let $\psi\in\Psi$ satisfy condition \eqref{good upper bound}. If $0\leq x\in M_{\psi},$ then there exists a constant $c(x)\in\mathbb{N}$ such that
$$n_x(\frac{\psi(t)}{t})\leq c(x)t$$
for every sufficiently large $t.$
\end{lem}
\begin{proof} Set $\varphi(t)=t/\psi(t).$ It follows from \eqref{good upper bound} that there exists a constant $\alpha>0$ and $t_0>0$ such that
$$\varphi(2t)\geq 2^{\alpha}\varphi(t)$$
for every $t\geq t_0.$ Thus,
$$\varphi(2^nt)\geq 2^{n\alpha}\varphi(t)$$
for $t\geq t_0.$

Consider sets $A$ and $B$ defined by the formula
$$A:=\{s:\ x^*(s)>\frac{\psi(t)}{t}\}\subset \{s: \|x\|_{M_{\psi}}\frac{\psi(s)}{s}>\frac{\psi(t)}{t}\}=:B.$$

Fix $c=2^n$ such that $2^{n\alpha}\geq\max\{1,\|x\|_{M_{\psi}}\}.$ It follows that $\varphi(ct)>\|x\|_{M_{\psi}}\varphi(t)$ for all $t\geq t_0.$ Therefore, $ct\notin B$ if $t\geq t_0.$ Since $\varphi$ is an increasing function (see \cite{KPS}), we have $\sup B\leq ct$ for $t\geq t_0.$ Since $B$ is an interval, we have $m(B)\leq ct$ provided that $t\geq t_0.$ Thus, for $t\geq t_0,$ we have $n_x(\psi(t)/t)=m(A)\leq m(B)\leq ct.$
\end{proof}

\begin{rem}\label{psi compatible remark} Let $\psi\in\Psi$ and let $\tau_{\omega}$ be a Dixmier trace on $M_{\psi}.$ We have
$$\omega(\frac{\psi(nt)}{\psi(t)})=1$$
for every $n\geq1.$ Indeed, if $\tau_{\omega}$ is linear then (see \cite{DPSS}) $$\omega(\frac{\psi(nt)}{\psi(t)})=\tau_{\omega}(n\sigma_{1/n}\psi')=\tau_{\omega}(\psi')=1.$$
\end{rem}

This remark is frequently used below together with the following lemma from \cite{DPSSS}.
\begin{lem}\label{34} Let $\omega\in L_{\infty}(0,\infty)^*$ be an arbitrary generalised limit. If $x,y\in L_{\infty}(0,\infty)$ are such that $\omega(|x-1|)=0,$ then $\omega(xy)=\omega(y).$
\end{lem}

\begin{lem}\label{ocenka snizu} Let $\psi\in\Psi$ satisfy the condition \eqref{good upper bound} and let $\tau_{\omega}$ be a Dixmier trace on $M_{\psi}.$ We have
$$\omega(\frac{-1}{\psi(t)}\int_{\psi(t)/t}^{\infty}\lambda dn_x(\lambda))\leq\tau_{\omega}(x)$$
for every positive $x\in M_{\psi}.$
\end{lem}
\begin{proof} Let $c(x)$ be the constant defined in Lemma \ref{dx estimate}. Clearly,
$$\frac{1}{\psi(t)}\int_0^{n_x(\psi(t)/t)}x^*(s)ds=(\frac{\psi(c(x)t)}{\psi(t)})\cdot(\frac{1}{\psi(c(x)t)}\int_0^{n_x(\psi(t)/t)}x^*(s)ds).$$
It follows from Remark \ref{psi compatible remark} and Lemma \ref{34} that
$$\omega(\frac{1}{\psi(t)}\int_0^{n_x(\psi(t)/t)}x^*(s)ds)=\omega(\frac{1}{\psi(c(x)t)}\int_0^{n_x(\psi(t)/t)}x^*(s)ds).$$
It follows from Lemma \ref{dx estimate} that
$$\omega(\frac{1}{\psi(c(x)t)}\int_0^{n_x(\psi(t)/t)}x^*(s)ds)\leq\omega(\frac{1}{\psi(c(x)t)}\int_0^{c(x)t}x^*(s)ds).$$
However, since $\omega$ is dilation invariant, we have
$$\omega(\frac{1}{\psi(c(x)t)}\int_0^{c(x)t}x^*(s)ds)=\omega(\frac{1}{\psi(t)}\int_0^{t}x^*(s)ds).$$
\end{proof}

\begin{lem}\label{ocenka sverhu} Let $\psi\in\Psi$ and let $\tau_{\omega}$ be a Dixmier trace on $M_{\psi}.$ We have
$$\tau_{\omega}(x)\leq\omega(\frac{-1}{\psi(t)}\int_{\psi(t)/t}^{\infty}\lambda dn_x(\lambda))$$
for every positive $x\in M_{\psi}.$
\end{lem}
\begin{proof} Fix $n\in\mathbb{N}.$ Clearly,
$$\frac1{\psi(t)}\int_0^tx^*(s)ds=(\frac{\psi(nt)}{\psi(t)})\cdot(\frac1{\psi(nt)}\int_0^tx^*(s)ds).$$
It follows from Remark \ref{psi compatible remark} and Lemma \ref{34} that
\begin{equation}\label{tauomega equality}
\tau_{\omega}(x)=\omega(\frac1{\psi(nt)}\int_0^tx^*(s)ds).
\end{equation}

We claim that
$$\int_0^tx^*(s)ds\leq\int_0^{n_x(\psi(nt)/nt)}x^*(s)ds+\frac1n\psi(nt).$$
The inequality is evident if $t\leq n_x(\psi(nt)/nt).$ If $t>n_x(\psi(nt)/nt),$ then $x^*(s)\leq\psi(nt)/nt$ for every $s\in[n_x(\psi(nt)/nt),t].$ Thus,
$$\int_0^tx^*(s)ds=\int_0^{n_x(\psi(nt)/nt)}x^*(s)ds+\int_{n_x(\psi(nt)/nt)}^tx^*(s)ds\leq$$
$$\leq\int_0^{n_x(\psi(nt)/nt)}x^*(s)ds+(t-n_x(\frac{\psi(nt)}{nt}))\cdot\frac{\psi(nt)}{nt}$$
and the claim follows.

Hence,
$$\omega(\frac1{\psi(nt)}\int_0^tx^*(s)ds)\leq\omega(\frac1{\psi(nt)}\int_0^{n_x(\psi(nt)/nt)}x^*(s)ds)+\frac1n.$$
It follows from \eqref{tauomega equality} and the dilation-invariance of $\omega$ that
$$\tau_{\omega}(x)\leq\omega(\frac1{\psi(t)}\int_0^{n_x(\psi(t)/t)}x^*(s)ds)+\frac1n.$$
Since $n$ is arbitrary large, we are done.
\end{proof}

The following theorem is the principal result of this section. It follows immediately from Lemmas \ref{ocenka snizu} and \ref{ocenka sverhu}.

\begin{thm}\label{lidskii comm} Let $\psi\in\Psi$ satisfy the condition \eqref{good upper bound} and let $\tau_{\omega}$ be a Dixmier trace on $M_{\psi}.$ We have
$$\tau_{\omega}(x)=\omega(\frac{-1}{\psi(t)}\int_{\psi(t)/t}^{\infty}\lambda dn_x(\lambda))$$
for every positive $x\in M_{\psi}.$
\end{thm}

Arguing similarly, one can obtain similar assertion for Marcinkiewicz sequence spaces.

\begin{thm}\label{lidskii comm seq} Let $\psi\in\Psi$ satisfy the condition \eqref{good upper bound} and let $\tau_{\omega}$ be a Dixmier trace on $m_{\psi}.$ We have
$$\tau_{\omega}(x)=\omega(\frac{1}{\psi(n)}\sum_{x_k\geq\psi(n)/n}x_k)$$
for every positive $x\in m_{\psi}.$
\end{thm}

\section{Adjusted Lidskii formula for Dixmier traces: noncommutative setting}\label{noncomm section}

In this section, we extend preceding results to Dixmier traces on Marcin\-kiewicz operator ideals.

\subsection{Adjusted Lidskii formula for Dixmier traces: normal operators}

The following assertion follows directly from the Theorem \ref{lidskii comm seq}.

\begin{lem}\label{self-adjoint} Let $\psi\in\Psi$ satisfy the condition \eqref{good upper bound} and let $\tau_{\omega}$ be a Dixmier trace on $\mathcal{M}_{\psi}.$ We have
$$\tau_\omega(S)=\omega(\frac1{\psi(n)}\sum_{\lambda\in\sigma(S),|\lambda|>\psi(n)/n}\lambda)$$
for every self-adjoint operator $S\in\mathcal{M}_{\psi}.$
\end{lem}

The following three lemmas are used to extend the formula above to the case of normal operators.

\begin{lem}\label{useful estimate} Let $\psi\in\Psi$ satisfy the condition \eqref{good upper bound} and let $\tau_{\omega}$ be a Dixmier trace on $M_{\psi}.$ We have
$$\omega(\frac1tn_x(\frac{\psi(t)}t))=0$$
for every positive $x\in M_{\psi}.$ A similar assertion holds for Marcinkiewicz sequence space $m_{\psi}.$
\end{lem}
\begin{proof} Fix $n\in\mathbb{N}.$ It follows from the dilation-invariance of $\omega$ that
\begin{equation}\label{another tauomega equality}
\omega(\frac1tn_x(\frac{\psi(t)}t))=\omega(\frac1{nt}n_x(\frac{\psi(nt)}{nt})).
\end{equation}
It is clear that
$$\frac1{nt}n_x(\frac{\psi(nt)}{nt})=\frac1n+\frac1{\psi(nt)}\int_t^{n_x(\psi(nt)/nt)}\frac{\psi(nt)}{nt}ds.$$
If $t>n_x(\psi(nt)/nt),$ we have
$$\int_t^{n_x(\psi(nt)/nt)}\frac{\psi(nt)}{nt}ds\leq0.$$
If $t\leq n_x(\psi(nt)/nt),$ then
$$\int_t^{n_x(\psi(nt)/nt)}\frac{\psi(nt)}{nt}ds\leq\int_t^{n_x(\psi(nt)/nt)}x^*(s)ds\leq\int_t^{c(x)nt}x^*(s)ds.$$
The last inequality holds for all sufficiently large $t$ by Lemma \ref{dx estimate}.

In either case,
$$0\leq\frac1{nt}n_x(\frac{\psi(nt)}{nt})\leq\frac1n+\frac1{\psi(nt)}\int_t^{c(x)nt}x^*(s)ds.$$
It follows now from the \eqref{another tauomega equality} that
$$\omega(\frac1tn_x(\frac{\psi(t)}t))\leq\frac1n+\omega(\frac1{\psi(nt)}\int_t^{c(x)nt}x^*(s)ds).$$
It is clear that
$$\omega(\frac1{\psi(nt)}\int_t^{c(x)nt}x^*(s)ds)=$$
$$=\omega(\frac1{\psi(nt)}\int_0^{c(x)nt}x^*(s)ds)-\omega(\frac1{\psi(nt)}\int_0^tx^*(s)ds).$$
It follows from the dilation-invariance of $\omega$ that
$$\omega(\frac1{\psi(nt)}\int_t^{c(x)nt}x^*(s)ds)=\omega(\frac1{\psi(t)}\int_0^{c(x)t}x^*(s)ds)-\omega(\frac1{\psi(nt)}\int_0^tx^*(s)ds).$$
It follows from Remark \ref{psi compatible remark} and Lemma \ref{34} that both terms in the right-hand side of the equality above are equal to $\tau_{\omega}(x).$

Therefore,
$$\omega(\frac1tn_x(\frac{\psi(t)}t))\leq\frac1n.$$
Since $n$ is arbitrary large, we are done.
\end{proof}

\begin{lem}\label{nenujnye chleny} Let $\psi\in\Psi$ satisfy the condition \eqref{good upper bound} and let $\tau_{\omega}$ be a Dixmier trace on $\mathcal{M}_{\psi}.$ We have
$$\omega(\frac1{\psi(n)}\sum_{|\Re\lambda|>\psi(n)/n,|\Im\lambda|\leq\psi(n)/n,\lambda\in\sigma(S)}\Im\lambda)=0,$$
$$\omega(\frac1{\psi(n)}\sum_{|\Re\lambda|\leq\psi(n)/n,|\Im\lambda|>\psi(n)/n,\lambda\in\sigma(S)}\Re\lambda)=0$$
for any normal operator $S\in\mathcal{M}_{\psi}.$
\end{lem}
\begin{proof} We prove the first assertion only. Proof of the second one is identical.

Note that $\lambda\in\sigma(S)$ if and only if $|\lambda|\in\sigma(|S|).$ It follows immediately that
$$|\sum_{|\Re\lambda|>\psi(n)/n,|\Im\lambda|\leq\psi(n)/n,\lambda\in\sigma(S)}\Im\lambda|\leq\sum_{|\Re\lambda|>\psi(n)/n,|\Im\lambda|\leq\psi(n)/n,\lambda\in\sigma(S)}\frac{\psi(n)}{n}\leq$$
$$\leq\frac{\psi(n)}{n}\sum_{|\lambda|>\psi(n)/n,\lambda\in\sigma(S)}1=\frac{\psi(n)}{n}\sum_{\lambda>\psi(n)/n,\lambda\in\sigma(|S|)}1=\frac{\psi(n)}{n}n_{|S|}(\frac{\psi(n)}{n}).$$
The assertion follows now from Lemma \ref{useful estimate}.
\end{proof}

\begin{lem}\label{eshe nenujnye chleny} Let $\psi\in\Psi$ satisfy the condition \eqref{good upper bound} and let $\tau_{\omega}$ be a Dixmier trace on $\mathcal{M}_{\psi}.$ We have
$$\omega(\frac1{\psi(n)}\sum_{|\Re\lambda|,|\Im\lambda|\leq\psi(n)/n,|\lambda|>\psi(n)/n,\lambda\in\sigma(S)}\lambda)=0$$
for any normal operator $S\in\mathcal{M}_{\psi}.$
\end{lem}
\begin{proof} It is clear that
$$|\sum_{|\Re\lambda|,|\Im\lambda|\leq\psi(n)/n,|\lambda|>\psi(n)/n,\lambda\in\sigma(S)}\lambda|\leq\sum_{|\Re\lambda|,|\Im\lambda|\leq\psi(n)/n,|\lambda|>\psi(n)/n,\lambda\in\sigma(S)}|\lambda|\leq$$
$$\leq\sum_{\psi(n)/n<|\lambda|\leq2\psi(n)/n,\lambda\in\sigma(S)}|\lambda|\leq\frac{2\psi(n)}{n}\sum_{|\lambda|>\psi(n)/n,\lambda\in\sigma(S)}1=$$
$$=\frac{2\psi(n)}{n}\sum_{\lambda>\psi(n)/n,\lambda\in\sigma(|S|)}1=\frac{2\psi(n)}{n}n_{|S|}(\frac{\psi(n)}{n}).$$
The assertion follows now from Lemma \ref{useful estimate}.
\end{proof}

The following theorem extends result of Lemma \ref{self-adjoint} to normal operators from $M_{\psi}.$

\begin{thm}\label{lidskiy normal} Let $\psi\in\Psi$ satisfy the condition \eqref{good upper bound} and let $\tau_{\omega}$ be a Dixmier trace on $\mathcal{M}_{\psi}.$ We have
$$\tau_\omega(S)=\omega(\frac1{\psi(n)}\sum_{|\lambda|>\psi(n)/n,\lambda\in\sigma(S)}\lambda)$$
for any normal operator $S\in\mathcal{M}_{\psi}.$
\end{thm}
\begin{proof} It follows from Lemma \ref{self-adjoint} that
$$\tau_{\omega}(\Re S)=\omega(\frac1{\psi(n)}\sum_{|\lambda|>\psi(n)/n),\lambda\in\sigma(\Re S)}\lambda)=\omega(\frac1{\psi(n)}\sum_{|\Re\lambda|>\psi(n)/n,\lambda\in\sigma(S)}\Re\lambda).$$
By Lemma \ref{nenujnye chleny},
$$\tau_{\omega}(\Re S)=\omega(\frac1{\psi(n)}\sum_{\max\{|\Re\lambda|,|\Im\lambda|\}>\psi(n)/n,\lambda\in\sigma(S)}\Re\lambda).$$
The same is valid for $\Im S.$ By the linearity,
$$\tau_{\omega}(S)=\omega(\frac1{\psi(n)}\sum_{\max\{|\Re\lambda|,|\Im\lambda|\}>\psi(n)/n,\lambda\in\sigma(S)}\lambda).$$
It follows from Lemma \ref{eshe nenujnye chleny} that
$$\tau_{\omega}(S)=\omega(\frac1{\psi(n)}\sum_{|\lambda|>\psi(n)/n,\lambda\in\sigma(S)}\lambda).$$
\end{proof}

\subsection{Adjusted Lidskii formula for Dixmier traces: general case}

Recall the following result of Ringrose (see Theorems 1, 6 and 7 from \cite{Ringrose}).

\begin{thm} Let $T\in B(H)$ be a compact operator. There exists a projection-valued measure $E_{\lambda}$ such that
\begin{enumerate}
\item $$TE_{\lambda}=E_{\lambda}TE_{\lambda}.$$
\item Either $E_{\lambda}=E_{\lambda-0}$ or
$${\rm rank}(E_{\lambda}-E_{\lambda-0})=1.$$
\item\label{quasinilpotence criterion} If, in addition,
$$TE_{\lambda}=E_{\lambda-0}TE_{\lambda},$$
then $T$ is quasi-nilpotent.
\end{enumerate}
\end{thm}

\begin{cor}\label{normal plus quasinilpotent} Let $T\in B(H)$ be a compact operator. There exist compact normal operator $S$ and compact quasi-nilpotent operator $Q$ such that $T=S+Q$ and $\sigma(S)=\sigma(T).$
\end{cor}

\begin{proof} Define an operator $S$ by the following formula
$$S=\sum_{E_{\lambda}\neq E_{\lambda-0}}(E_{\lambda}-E_{\lambda-0})T(E_{\lambda}-E_{\lambda-0}).$$
A straightforward computation shows that the operator $Q=T-S$ satisfies the condition \ref{quasinilpotence criterion} of the Theorem above. Hence, $Q$ is quasi-nilpotent.

Evidently, $S$ is a diagonal operator with eigenvalues of $T$ on the diagonal. Hence, $\sigma(S)=\sigma(T).$
\end{proof}

By the Weil theorem, sequence of eigenvalues of $T$ is majorized by the sequence of its singular values (see Theorem 3.1 from \cite{GohbergKrein}). Hence, for $T\in \mathcal{M}_{\psi},$ we obtain $S,Q\in\mathcal{M}_{\psi}.$

The following assertion directly follows from the Theorem 3.3 from \cite{Kalton}).

\begin{thm} If $Q\in\mathcal{M}_{\psi}$ is a quasi-nilpotent operator, then $Q$ belongs to the commutator $[\mathcal{M}_{\psi},B(H)].$
\end{thm}

\begin{cor}\label{quasinilpotent inessential} If $Q\in\mathcal{M}_{\psi}$ is a quasi-nilpotent operator and $\tau_{\omega}$ is an arbitrary Dixmier trace on $M_{\psi},$ then $\tau_{\omega}(Q)=0.$
\end{cor}
Indeed, due to \cite{Connes}, we have $\tau_{\omega}([A,B])=0$ for every $A\in \mathcal{M}_{\psi}$ and every $B\in B(H).$

The following theorem is the main result of this section.

\begin{thm}\label{lidskii formula} Let $\psi\in\Psi$ satisfy condition \eqref{good upper bound} and let $\tau_{\omega}$ be a Dixmier trace on $\mathcal{M}_{\psi}.$ We have
$$\tau_\omega(T)=\omega(\frac1{\psi(n)}\sum_{\lambda\in\sigma(T),|\lambda|>\psi(n)/n}\lambda)$$
for any operator $T\in\mathcal{M}_{\psi}.$
\end{thm}
\begin{proof} Let $S$ be a normal operator constructed in Corollary \ref{normal plus quasinilpotent}. The assertion holds for $S$ by Theorem \ref{lidskiy normal}. Note that $\tau_{\omega}(T)=\tau_{\omega}(S)$ by Corollaries \ref{normal plus quasinilpotent} and \ref{quasinilpotent inessential}. Since $\sigma(S)=\sigma(T),$ we are done.
\end{proof}

\section{Applications to heat kernel formula}

In this section, we provide a simple proof of one of the heat semigroup formulae from \cite{Sedaev} (see also earlier results in \cite{CPS,CRSS}). Our hypothesis on $\omega$ is very mild.

\begin{lem}\label{another M lemma} For any positive $x\in M_{1,\infty}$ we have
$$M(\frac1{\log(t)}\int_t^{n_x(1/t)}(x^*(s)-1/t)ds)=o(1).$$
\end{lem}
\begin{proof} If $t>n_x(1/t),$ we have
$$|\int_t^{n_x(1/t)}(x^*(s)-1/t)ds|\leq 1.$$
If $t\leq n_x(1/t),$ then $x^*(s)\geq1/t$ for every $s\in[t,n_x(1/t)].$ Therefore,
$$0\leq\int_t^{n_x(1/t)}(x^*(s)-1/t)ds\leq\int_t^{n_x(1/t)}x^*(s)ds.$$
The assertion follows now from the Lemma \ref{M k nulyu}.
\end{proof}

\begin{thm} Let $\tau_{\omega}$ be a Dixmier trace on $\mathcal{M}_{1,\infty}$ such that $\omega=\omega\circ M.$ We have
$$\tau_{\omega}(T)=\frac{\alpha}{\Gamma(1/\alpha)}\omega(\frac1t\sum_{\lambda\in\sigma(T)}\exp(-(t\lambda)^{-\alpha}))$$
for every positive operator $T\in\mathcal{M}_{1,\infty}.$
\end{thm}
\begin{proof} Let $x=x^*\in M_{1,\infty}$ be the rearrangement of $T,$ that is $x=i(\{s_n(T)\}).$ Without loss of generality, $x\leq 1.$ Since distributions of $T$ and $x$ coincide, we have
$$\omega(\frac1t\sum_{\lambda\in\sigma(T)}\exp(-(t\lambda)^{-\alpha}))=\omega(\frac1t\int_0^{\infty}\exp(-(tx(s))^{-\alpha})ds).$$
Setting $1/x(s)=u,$ we obtain
$$\omega(\frac1t\int_0^{\infty}\exp(-(tx(s))^{-\alpha})ds)=\omega(\frac1t\int_0^{\infty}\exp(-(u/t)^{\alpha})dn_x(1/u)).$$
It follows from the weak version of Karamata Theorem (see \cite{CPS,Sedaev} for details) that
$$\frac{\alpha}{\Gamma(1/\alpha)}\omega(\frac1t\sum_{\lambda\in\sigma(T)}\exp(-(t\lambda)^{-\alpha}))=\omega(\frac1tn_x(1/t)).$$
It is clear that
$$M^2(\frac1tn_x(1/t))-M(\frac1{\log(t)}\int_1^tx(s)ds)=$$
$$=M(\frac1{\log(t)}(\int_1^t\frac1{s^2}n_x(1/s)ds-\int_1^tx(s)ds)).$$
Integrating by parts, we obtain
$$\int_1^t\frac1{s^2}n_x(1/s)ds=-\frac1tn_x(1/t)+\int_1^t\frac1sdn_x(1/s)=\int_1^{n_x(1/t)}x(s)ds-\frac1tn_x(1/t).$$
Hence,
$$\int_1^t\frac1{s^2}n_x(1/s)ds-\int_1^tx(s)ds=\int_t^{n_x(1/t)}(x(s)-1/t)ds-1.$$
It follows from the Lemma \ref{another M lemma} that
$$M^2(\frac1tn_x(1/t))-M(\frac1{\log(t)}\int_1^tx(s)ds)=o(1).$$
The assertion follows now from the $M-$invariance of $\omega.$
\end{proof}


\begin{thebibliography}{100}
\bibitem{AS} Azamov N., Sukochev F. {\it A Lidskii type formula for Dixmier traces}, C.R. Math. Acad. Sci. Paris 340 (2005), no. 2, 107--112.
\bibitem{BF} Benameur M., Fack T. {\it Type II non-commutative geometry. I. Dixmier trace in von Neumann algebras},  Adv. Math.  199  (2006),  no. 1, 29--87.
\bibitem{CPS} Carey A., Phillips J., Sukochev F. {\it Spectral flow and Dixmier traces}, Adv. Math. 173 (2003), no. 1., 68--113.
\bibitem{CPRS1} Carey A., Phillips J., Rennie A., Sukochev F. {\it The Hochschild class of the Chern character for semifinite spectral triples},  J. Funct. Anal., vol.~213, (2004) no.~1, ~111--153.
\bibitem{CRSS} Carey A., Rennie A., Sedaev A., Sukochev F. {\it The Dixmier trace and asymptotics of zeta functions}, J. Funct. Anal. 249 (2007), no. 2, 253--283.
\bibitem{CS} Carey A., Sukochev F. {\it Dixmier traces and some applications in non-commutative geometry}, Russian Math. Surveys 61:6 1039--1099.
\bibitem{Connes} Connes A. {\it Noncommutative Geometry}, Academic Press, San Diego, 1994.
\bibitem{Dixmier} Dixmier J. {\it Existence de traces non normales}, C. R. Acad. Sci. Paris, 262 (1966).
\bibitem{DPSSS} Dodds P., de Pagter B., Sedaev A., Semenov E., Sukochev F. {\it Singular symmetric functionals}, (Russian)  Zap. Nauchn. Sem. S.-Peterburg. Otdel. Mat. Inst. Steklov. (POMI)  290  (2002),  Issled. po Linein. Oper. i Teor. Funkts. 30, 42--71, 178;  translation in  J. Math. Sci. (N. Y.)  124  (2004),  no. 2, 4867--4885
\bibitem{DPSSS2} Dodds P., de Pagter B., Sedaev A., Semenov E., Sukochev F. {\it Singular symmetric functionals and Banach limits with additional invariance properties}, (Russian) {\em Izv. Ross. Akad. Nauk Ser. Mat.}, vol.~67, no.~6 (2003) ~111--136.
\bibitem{DPSS} Dodds P., de Pagter B., Semenov E., Sukochev F. {\it Symmetric functionals and singular traces}, Positivity  2  (1998),  no. 1, 47--75.
\bibitem{GohbergKrein} Gohberg I., Krein M. {\it Introduction to the theory of linear nonselfadjoint operators}, Translations of Mathematical Monographs, Vol. 18 American Mathematical Society, Providence, R.I. 1969
\bibitem{Kalton} Kalton N. {\it Spectral characterization of sums of commutators. I}, J. Reine Angew. Math. 504 (1998), 115--125.
\bibitem{KPS} Krein S., Petunin Ju. and Semenov E. {\it Interpolation of linear operators}, Nauka, Moscow, 1978 (in Russian); English translation in Translations of Math. Monographs, Vol. {\bf 54}, Amer. Math. Soc., Providence, RI, 1982.
\bibitem{LSS} Lord S., Sedaev A., Sukochev F. {\it Dixmier traces as singular symmetric functionals and applications to measurable operators}, J. Funct. Anal. 224 (2005), no. 1, 72--106.
\bibitem{Pie} Pietsch A. {\it About the Banach Envelope of $l_{1,\infty}$}, Rev. Mat. Complut. 22~(1) (2009) 209--226.
\bibitem{Ringrose} Ringrose J. {\it Super-diagonal forms for compact linear operators},  Proc. London Math. Soc. (3) 12 (1962) 367--384.
\bibitem{Sedaev} Sedaev A. {\it Generalized limits and related asymptotic formulas}, Math.Notes. 86:4 (2009), 612-627.
\end{thebibliography}
\end{document}